\documentclass [a4paper,12pt] {article}
\usepackage [latin1]{inputenc}
\usepackage[all]{xy}
\usepackage [english]{babel}
\usepackage{amsmath}
\usepackage{amsfonts}
\usepackage{amssymb}
\usepackage{graphics}
\usepackage{graphicx}
\newtheorem{teorema}{Theorem}[section]
\newtheorem{lemma}[teorema]{Lemma}
\newtheorem{propos}[teorema]{Proposition}
\newtheorem{corol}[teorema]{Corollary}
\newtheorem{ex}{Example}[section]
\newtheorem{rem}{Remark}[section]
\newtheorem{defin}[teorema]{Definition}
\def\bt{\begin{teorema}}
\def\et{\end{teorema}}
\def\bp{\begin{propos}}
\def\ep{\end{propos}}
\def\bl{\begin{lemma}}
\def\el{\end{lemma}}
\def\bc{\begin{corol}}
\def\ec{\end{corol}}
\def\br{\begin{rem}\rm}
\def\er{\end{rem}}
\def\bex{\begin{ex}\rm}
\def\eex{\end{ex}}
\def\bd{\begin{defin}}
\def\ed{\end{defin}}
\def\demo{\par\noindent{\bf Proof.\ }}
\def\enddemo{\ $\Box$\par\vskip.6truecm}
\def\R{{\mathbb R}}   \def\a {\alpha} \def\g{\gamma}
\def\N{{\mathbb N}}     \def\d {\delta} \def\e{\varepsilon}
\def\C{{\mathbb C}}      \def\l{\lambda}
                 
                 \def\t{\theta}  \def\z{\zeta}
\def\P{{\mathbb P}}
\def\Re{{\sf Re}}

 \def\oli{\overline}
\def\D{\Delta}
\def\G{\Gamma}
\def\L{\Lambda}
\def\O{\Omega}
\def\S{\Sigma}
\begin{document}
\title {Liouville-type theorems for foliations with complex leaves}
\author {Giuseppe Della Sala}
\maketitle

\section*{Introduction}

In this paper we discuss some structure problems regarding the
foliation of Levi flat manifolds, i.e.\ those $CR$ manifolds which
are foliated by complex leaves or, equivalently, whose Levi form
vanishes. Levi flat manifolds have a particular significance among
$CR$ manifolds, since - due to their degenerate nature - they often
behave as \lq\lq limit case\rq\rq, or they are an obstruction in
extension problems of $CR$ objects (see, for example, \cite{LT}).
The geometry of Levi flat manifolds is a source of many interesting
problems (see \cite{Bry}, \cite{BI}).

Here, we address the following general question: assuming that a
Levi flat submanifold $S\subset \C^n$ is bounded in some directions,
what can be said about its foliation? We will show that, in some
circumstances, it is possible to conclude that the foliation (or
even, more generally, a complex leaf of a foliated manifold) is
\lq\lq trivial\rq\rq, i.e. is made up by complex planes. A first
result in this direction is the following \bt \label{coroll} Let $S$
be a smooth Levi-flat hypersurface of $\C^n=\C^{n-1}\times \C_w$,
contained in $C=\{|w|<1\}$ and closed in $\oli C$. Then $S$ is
foliated by hyperplanes $\{w=const.\}$.\et In order to treat this
problem it is useful to consider $S$ as an \emph{analytic
multifunction}. These objects, which were first introduced by Oka
\cite{Oka}, are set-valued functions $\C\to \mathbf k(\C)$ (where
$\mathbf k(\C)$ denote the subset of the power set $\mathcal P(\C)$
formed by the compact subsets of $\C$) which behave in some ways as
analytic function; namely, according to Oka's definition, the
complementary of their graph is pseudoconvex. However, we will find
more convenient to adopt as a definition the characterization found
by Slodkowski \cite{Slo} by means of plurisubharmonic (psh)
functions.

We shall see that Theorem \ref{coroll} becomes then a rather easy
consequence of the results already obtained for analytic
multifunctions, namely, the extension of Liouville's theorem to such
objects. Then, we generalize Theorem \ref{coroll} to higher
codimension (see Theorem \ref{coroll2}) by a slightly less trivial
application of the Liouville Theorem.

Later on, we will discuss other related problems which cannot be
treated by means of analytic multifunctions. In section
\ref{folgraf}, we consider the case of a (complex leaf of a)
foliation of the graph of a bounded function on $\C^n\times \R$. In
this case, an analysis of each single complex leaf is required.
Afterwards, in section \ref{folcil} we consider a (real) foliation
of $D\times \C$ by complex leaves, and we show that under suitable
(maybe too restrictive) geometric conditions on this foliation it is
again possible to prove a triviality result.

\section*{Acknowledgements}
This paper would not have been written without the help of
G.Tomassini, who suggested the problems and gave many important
advices.

\noindent I also wish to thank N.Shcherbina for several helpful
talks, expecially for the idea of adopting the methods of analytic
multifunctions.

\section{Levi flat manifolds contained in a cylinder} \label{uno}

\subsection{Analytic multifunctions and Liouville Theorem} Consider a
function $f:\C^n\to \mathcal P(\C)$, i.e.\  a set-valued function
from $\C^n$ to the power set of $\C$. Let $\G(f)\subset \C^{n+1}$ be
defined as
$$\G(f)=\bigcup_{z\in\C^n} \{z\}\times f(z).$$
We say that $f$ is an \emph{analytic multifunction} if each value
$f(z)$ is a compact set and $\C^{n+1}\setminus \G(f)$ is
pseudoconvex. With this definition, a holomorphic function
$f\in\mathcal O(\C^n)$ is clearly an analytic multifunction.

Let $\rho:\C^{n+1}\to \R$ be a continuous plurisubharmonic function.
Let $\rho':\C^n\to \R$ be defined as
$$\rho'(z)=\max_{w\in f(z)} \rho(w).$$ In \cite{Slo} the following is proved
: \bl \label{def2} For any analytic multifunction $f$ and continuous
psh function $\rho$, $\rho'$ is a plurisubharmonic function.\el From
now on, by analytic multifunction we mean a multifunction for which
the conclusion of Lemma \ref{def2} holds true.

The following Liouville result (see also \cite{Ran}) on analytic
multifunction depends only on the property of Lemma \ref{def2}:

\bl \label{lioumul} Let $f$ be an analytic multifunction on $\C^n$,
and suppose that $f$ is bounded in the following sense:
$$\G(f)\subset\{|w|<M\}\subset\C^{n+1} $$
for some $M>0$. Let $\hat f$ be the multifunction defined as
$$\hat f(z)=\widehat{f(z)},\ z\in\C^n$$
where $\widehat K$ is the polynomial hull of $K$. Then $\hat f$ is
constant. \el \demo Let $P(w)$ be a polynomial on $\C_w$, and denote
again by $P$ the trivial extension of $P$ to $\C^{n+1}$
$P(z,w)=P(z)$. Then $|P|$ is a plurisubharmonic function on
$\C^{n+1}$, therefore by Lemma \ref{def2}
$$P'(z)=\max \{|P(w)|: w\in f(z)\}$$
is psh on $\C^n$. But, defining
$$C=\max_{|w|\leq M} P(w) $$
we have that $P'(z)\leq C$ for all $z\in \C^n$. Then, by Liouville's
Theorem for psh functions it follows that $P'$ is constant. We
deduce that $\hat f$ is constant. Indeed, in the opposite case we
could find $w_1\in \C$ and $z_1, z_2\in \C^n$ such that $w_1\in
(\hat f(z_1)\setminus \hat f(z_2))$, i.e. there would exist a
polynomial $P_1$ such that $$|P_1(w_1)|>\max_{\hat f(z_2)}|P_1|,$$
hence $$P_1'(z_2)<|P_1(w_1)|\leq P_1'(z_2)$$ which is a
contradiction.
\enddemo

\bex The hypothesis of Lemma \ref{lioumul} does not imply that $f$
is in turn a constant multifunction. A simple example is the
following:
$$f(z)=\left\{
    \begin{array}{ll}
      \{|w|=1\}, & \hbox{$z\neq 0$;} \\
      \{|w|\leq 1\}, & \hbox{$z=0$.}
    \end{array}
  \right.
 $$
\eex \bex A modification of the previous example shows that, even if
$\G(f)$ is a (disconnected) manifold, $f$ need not be constant if we
adopt the second definition of analytic multifunction (i.e. the
property discussed in Lemma \ref{def2}). Indeed, define $f(z)$ to be
the union of the unit circle $bD$ and any compact set contained in
the unit disc $D$; then, since any subharmonic function can \lq\lq
detect\rq\rq\ the behavior of $f$ only in $bD$, $f$ satisfies the
statement of Lemma \ref{def2} but the complement is not
pseudoconvex. As we show below, nevertheless, the result holds if
$\G(f)$ has the structure of a (even disconnected) Levi flat
manifold (which is obviously not the case in the previous example).
\eex

Lemma \ref{lioumul} provides a useful tool which allows to prove
Theorem \ref{coroll} quite easily. In fact, setting
$$f_S(\zeta)=S\cap \{z=\zeta\}$$
for $\zeta\in \C^n$ we have that $f_S$ is by definition an analytic
multifunction.

\noindent \textbf{Proof of Theorem \ref{coroll}}. By hypothesis the
multifunction $f_S$ is bounded, therefore in view of Liouville's
Theorem $\hat f_S$ is constant. We have to show that the
multifunction $f_S$ is constant, too. In order to do this, choose
$z_0\in \C^n$ in such a way that the complex line
$L_{z_0}=\{z=z_0\}\subset \C^{n+1}$ intersects $S$ transversally.
This means that $f(z_0)$ is a smooth compact real $1$-submanifold of
$\C$, i.e. a finite set $\{\l_i(z_0)\}_{1\leq i\leq k(z_0)}$ of
simple $C^\infty$ loops contained in $D=\{|w|<1\}$. Let $U_i(z_0)$
be the bounded connected component of $\C\setminus \l_i(z_0)$, and
let $\{\a_j(z_0)\}_{1\leq j \leq h(z_0)}$ be the \lq\lq
maximal\rq\rq loops, i.e. those $\l_i$'s which are not contained in
any $U_j$. For every $z\in \C^n$ such that $L_z\cap S$ is
transversal, adopting the same notations as above we define
$$\mathcal M(z)=\bigcup_{1\leq i\leq h(z)} \a_i(z).$$
Let $I\subset \C^n$ be the set of $z\in\C^n$ for which
\begin{itemize}
  \item $L_z$ has transversal intersection with $\bigcup_\zeta \mathcal
M(\zeta)$;
  \item $\mathcal M(z) = \mathcal M(z_0)$.
\end{itemize}
It suffices to show that $I$ is both open and closed. Indeed, in
this case $f'_S(z)= f_S(z)\setminus \mathcal M(z)$ is an analytic
multifunction, thus we can prove the statement of \ref{coroll}
inductively, where the induction is performed on the number of loops
of $f(z_0)$.

\underline{$I$ is open}. Let $z_1\in I$; clearly there exists a
neighborhood $\O$ of $z_1$ such that $h(z)=h(z_1)\equiv h$ for $z\in
\O$ and $\mathcal M_\O= \bigcup_{z\in\O} \mathcal M(z)$ is a
submanifold of $\O\times \C_w$ for which $L_z\cap \mathcal M_\O$ is
transversal. Moreover, observe that if $\{V_i(z)\}_{1\leq i\leq h}$
are the connected components of $\C\setminus \a_i(z)$, then $\hat
f(z) = \bigcup_i V_i(z)$. This implies immediately that $\mathcal
M(z)$ is constant on $\O$.

\underline{$I$ is closed}. Let $\mathcal M_I= \bigcup_{z\in I}
\mathcal M_z$ and let $z_2\in \oli I$. Then, $\oli {\mathcal
M_I}\cap L_{z_2}=\mathcal M(z_0)$; moreover, since $S$ is a smooth
manifold, we have
$$T_{(z_2,w')}(S)\supset \{(z,w)\in \C^{n+1}: w=w'\}$$
for every $w'\in f(z_0)$. But, since we clearly have $\oli {\mathcal
M_I}\cap L_{z_2}=\mathcal M(z_2)$, this implies that $L_{z_2}\cap
\mathcal M$ is transversal, i.e. $z_2\in I$.
\enddemo

\subsection{Higher codimension} The analogous of Theorem \ref{coroll}
for Levi flat surfaces $S$ of higher codimension can also be proved
by following the same methods. However, in this case, analyticity of
the multifunction $f_S$ defined by $S$ is more involved, due to the
fact that $S$ is no longer pseudoconvex. Also the proof of the fact
that $f_S$ is constant whenever $\hat f_S$ is needs to be adapted
using \cite{Sto} (see the proof of Theorem \ref{coroll2}).

We consider a real ($2d-1$)-codimensional submanifold $S\subset
\C^{n+d}\cong \C^n \times \C^d$, with coordinates
$z_1,\ldots,z_n,w_1,\ldots,w_d$.
 \bt \label{coroll2} Let $S\subset \C^{n+d}$ be a ($2d-1$)-codimensional
Levi flat submanifold (i.e. foliated by complex leaves of dimension
$n$), contained in
$$C=\{(z,w)\in \C^{n+d}:\sum_{i=1}^d |w_i|^2 < 1\}$$ and closed in
$\oli C$. Then $S$ is foliated by coordinate complex $n$-planes
$\{w_1=c_1, \ldots, w_d=c_d\}$. \et Let $f:\C^n\to \mathcal P(\C^d)$
be a function from $\C^n$ to the subsets of $\C^d$, $d\geq 2$. We
recall that, according to our definition, $f$ is an analytic
multifunction if $f(z)$ is compact for each $z\in \C^n$ and, for
every continuous plurisubharmonic function $\rho:\C^{n+d}\to \R$,
the function $\rho':\C^n\to \R$ defined as
$$\rho'(z)=\max_{f(z)}\rho(z,w)$$ is plurisubharmonic (see Lemma \ref{def2}).

Let $L_z$, $z\in \C^n$, be the vertical complex $d$-plane over $z$
i.e.
$$L_z=\{(\zeta,w)\in \C^{n+d}: \zeta=z\}.$$ Consider the
set-valued function $f_S$ defined by $f_S(z)=L_z\cap S$
(generically, $f_S(z)$ is the union of a finite number of loops). We
want to show that $f_S$ is an analytic multifunction.

\bl \label{fSmul}$f_S$ is an analytic multifunction.\el \demo Let
$\rho:\C^{n+d}\to \R$ be a psh function, and define $\rho'$ as
above. Let $z_0\in \C^n$, and let $\mathcal L \subset \C^n$ be a
complex line passing through $z_0$. For a generic choice of
$\mathcal L$, the intersection of $S$ with the complex $(d+1)$-plane
$$\{(z,w)\in \C^{n+d}:z\in\mathcal L\}$$
is transversal, and thus a Levi flat submanifold of $\C^{d+1}$.
Therefore, since it is sufficient to show that the restriction of
$\rho'$ to a generic $\mathcal L$ is subharmonic, we can suppose
$n=1$.

Assume, then, that $f_S$ is a $\mathcal P(\C^d)$-valued
multifunction defined over $\C_z$, and fix $z_0\in \C$. If $w\in
f(z_0)$, we denote by $\S_w$ the leaf of the foliation of $S$
through $w$. Two cases are possible: \

\

\begin{itemize}
  \item $T_{(z_0,w)}(\S_w)\nsubseteq \C^d$;
  \item $T_{(z_0,w)}(\S_w)\subset \C^d$.
\end{itemize}
  In the former, for a
  sufficiently small neighborhood $V_w=(\D\times U)_w$ of $(z_0,w)$ we have
  that $\S_w\cap V_w$ can be written
  $$\S_w\cap V_w =\{(z,w)\in \D\times U: w_1=g_1^w(z),\ldots,w_d=g_d^w(z) \}$$
  for some holomorphic function $g_i^w\in\mathcal O(\D)$. Moreover,
  observe that for $w'\in f(z_0)$ in a small enough neighborhood $\mathcal W_w$ of
  $w$, we can choose a $\D$ which does not depend on $w'$.

  In the latter, consider the
  restriction of the projection $\pi:\C^{d+1}\to \C$ to a small
  neighborhood $\mathcal V_w$ of $(z_0,w)$ in $\S_w$. We can suppose
  that $\mathcal V_w$ is a local chart such that $(z_0,w)\cong 0$. Denote by $\zeta$
  the complex coordinate on $\mathcal V_w$. Since $\pi|_{\mathcal V_w}$ is a holomorphic
  function and its prime derivative vanishes
  at $0$,  there exists $k\geq 1$ such that
  $$\frac{\partial^k}{\partial \zeta}\pi|_{\mathcal V_w}=0, \  \frac{\partial^{k+1}}{\partial \zeta}\pi|_{\mathcal V_w}\neq
  0.$$ Otherwise, we would have $\pi|_{\mathcal V_w}\equiv z_0$ and
  thus $\S_w$ would be a complex line contained in $\C^d$, which is
  impossible since it must be contained in the cylinder $C$ of Theorem
  \ref{coroll2}. It follows that $\pi|_{\mathcal V_w}$ is a
  $k$-sheeted covering over some neighborhood $\D$ of $z_0$. Now,
  the restriction of $\pi$ to the leaves $\S_{w'}$ passing through
  the points $(z_0,w')$ of a small neighborhood of $(z_0,w)$ can be
  interpreted as a smooth one-parameter family of holomorphic functions
  $\pi_t:\mathcal V_k\to \C_z$, such that $\pi_0=\pi$. For $|t|<<1$,
  the argument principle implies that the sum of the orders of the zeroes of $(\partial /\partial
  \zeta)\pi_t$ is still $k-1$. This in turn means that for $w'$
  sufficiently close to $w$ the projection $\pi|_{\S_{w'}}$ is still a
  $k$-sheeted covering over some neighborhood $\D_{w'}$; in a
  sufficiently small neighborhood $\mathcal W_w$ we can
  assume to have chosen a $\D$
  independent from $w'$.

  \

Since $f(z_0)$ is a compact set, we can choose finitely many open
sets as above, $\mathcal W_{w_1},\ldots,\mathcal W_{w_h}$, in such a
way that
$$\bigcup_{i=1}^h \mathcal W_{w_i} = f(z_0).$$
Choose a disc $\D\subset \D_{w_1}\cap \ldots \cap \D_{w_h}$. We
claim that $\rho'$ is plurisubharmonic on $\D$. In order to prove
this, choose $w\in f(z_0)$:
\begin{itemize}
  \item if $w\in \mathcal W_{w_j}$ with $w_j$ of the first kind,
  then we define
  $$\rho_w^j=\rho|_{\S_w\cap \pi^{-1}(\D)};$$
  \item if $w\in \mathcal W_{w_j}$ with $w_j$ of the second kind,
  we define
  $$\rho_w^j=(\max_{\S_w\cap \pi^{-1}(\D_{w_j})}\rho(z,w))|_\D.$$
\end{itemize}
In both cases, $\rho_w^j$ is a psh function. Observe that it may
occur $\rho_w^i\neq \rho_w^j$ when $i\neq j$. Nevertheless, consider
\begin{equation} \label{max} \varrho(z)= \max_{1\leq i\leq h, w\in f(z_0)}
\rho_w^i; \end{equation} we have to show that $\varrho(z) =
\rho'(z)$. A priori, the maximum in equation (\ref{max}) may be
performed, for $z\neq z_0$, on a proper subset of $f(z)$, due to the
possible existence of leaves of $S$ which accumulate to $f(z_0)$
without intersecting it. However, this does not happen. This is a
consequence of the fact that
$$\bigcup_{w\in f(z_0)}
\S_w\cap \pi^{-1}(\D)= S\cap \pi^{-1}(\D)$$ as proved in Lemma
\ref{sotto} below. It follows that $\varrho(z) = \rho'(z)$. Since we
already know that $\rho'(z)$ is continuous, (\ref{max}) implies
$\rho'(z)$ is plurisubharmonic.
\enddemo

\bl \label{sotto} Let $z_0\in \C_z$, and let $\S$ be a leaf of the
foliation of $S$ such that $z_0\in \oli{\pi(\S)}$. Then $z_0\in
\pi(\S)$. \el \demo Observe that, since $\S\subset C$, and
consequently the $w$-coordinates are bounded on $\S$, there exists
$w\in f(z_0)$ such that $(z_0,w)$ is a cluster point for $\S$. Take
a neighborhood $\mathcal U$ of $(z_0,w)$ such that the foliation is
trivial on $S\cap \mathcal U$. Then $\S\cap \mathcal U$ is a union
of leaves of this trivial foliation. Let $\S_0$ be the leaf of
$S\cap \mathcal U$ passing through $(z_0,w)$; then, either $\S$
contains $\S_0$ or it contains a sequence of leaves that converges
to $\S_0$. Suppose that the second case occurs. Arguing as in the
previous Lemma, we show that $\S_0$ is not \lq\lq vertical\rq\rq\
i.e. it is not contained in $\C^d$. Then $\pi(\S_0)$ is a open set
containing $z_0$; but for $\S'$ close enough to $\S_0$, $\pi(\S')$
is an open set containing $z_0$. This proves the thesis. \enddemo

Lemma \ref{fSmul} allows to prove, exactly in the same way as
before, that $\hat f_S$ is a constant multifunction. We have to
show, again, that this fact forces $f_S$ to be constant.

\noindent \textbf{Proof of Theorem \ref{coroll2}}. Observe that, for
$z$ belonging to a dense, open subset $J$ of $\C^n$, $L_z$
intersects $S$ transversally. For $z\in J$, $f(z)=L_z\cap S$ is the
disjoint union of a finite set $\{\g_i(z)\}_{1\leq i\leq k(z)}$ of
loops in $\C^d$. It is a well-known fact (\cite{Sto}) that, in this
case, the polynomial hull $\hat f(z)$ of $f(z)$ is given by the
union of some of the loops $\g_i$ and some complex varieties $\L_j$
whose boundaries are the others $\g_i$'s. We choose the minimal
subsets of loops $\{\a_i(z)\}_{1\leq i\leq h(z)}$ such that, if
$\mathcal M(z)= \a_1\cup\ldots\cup\a_{h(z)}$, then $\widehat
{\mathcal M}(z)=\hat f(z)$. Observe that $\mathcal M(z)$ is
univocally defined. It is sufficient to prove that $\mathcal M(z)$
is constant, because in such a case we can proceed inductively as in
the proof of Theorem \ref{coroll}. In view of the structure of the
hull $\hat f_S(z)$, it is clear that $\mathcal M(z)$ is constant for
$z\in J$. Moreover, arguing again as in \ref{coroll}, it is clear
that
$$\mathcal M=\bigcup_{z}\mathcal M(z)$$ is a manifold with
transversal intersection with $L_z$ also for $z\in \oli J$. It
follows that $\mathcal M$, and therefore $f_S$, is constant.
\enddemo

\section{Foliation of a graph} \label{folgraf}

Consider, in $\C^2$, coordinates $(z,w)$ with $z=x+iy$ and $w=u+iv$.
We denote by $\pi$ the projection $\pi:\C^2\to \C_z$ and by $\tau$
the projection $\tau:\C^2\to \C_z \times \R_u$. Let $\rho:\C\times
\R \to \R$ be a smooth function such that
$$S=\{v=\rho(z,u)\}$$
is a Levi flat graph. Then the leaves of the foliation associated to
$S$ are regular (immersed) complex manifolds of dimension $1$. We
claim that the following holds true: \bt \label{main7}Suppose that
$|\rho|$ is bounded by some constant $M$; then $S$ is foliated by
complex lines, i.e. $\rho$ has the form $$\rho (z,u) = \varrho(u)$$
for some smooth function $\varrho:\R\to\R$. \et Actually, we are
going to prove a sharper result. If $S$ is any foliated
$3$-submanifold, we say that a leaf $\S$ of $S$ is \emph{properly
embedded} if, for (almost) every ball $B\subset \C^2$, the connected
components of $\oli B\cap \S$ are compact, embedded submanifolds of
$\oli B\cap S$ with boundary. We say that the foliation of $S$ is
\emph{proper} if the leaves are properly embedded. \bt \label{main4}
Let $S=\{v=\rho(z,u)\}$ be a properly foliated hypersurface of
$\C^2$, and suppose that $\rho$ is bounded. Then every complex leaf
of $S$ is a complex line.\et Theorem \ref{main7} is a consequence of
Theorem \ref{main4}, the foliation of a Levi flat graph being proper
(see, for example, the main Theorem of \cite{Shc}). A simple proof
of \ref{main7} can be given using analytic multifunctions.

\noindent \textbf{Proof of Theorem \ref{main7}}. By hypothesis there
exists a complex line $\{w=c\}$ such that $S$ lies outside the
cylinder $$C=\{(z,w): |w-c|<\varepsilon\}.$$ Then, we can perform a
rational change of coordinates (acting only on the $w$-coordinate)
such that the image $S'$ of $S$ is contained in $\C_z\times D_w$,
where $D_w$ is the unit disc. The complement of $$\oli S'= S'\cup
(\C_z\times \{0\})$$ in $\C_z\times D_w$ is pseudoconvex. Indeed, a
psh exhaustion function $\varphi$ for $S$ on $\C^2$ induces a psh
exhaustion function $\varphi'$ for $S'$ on $\C_z\times
(D_w\setminus\{0\})$; then
$$\psi = \max \{\varphi' \left|\frac{1}{w}\right|\}$$
is a psh exhaustion function for $\oli S'$ on $\C_z\times D_w$. The
rest of the proof can be carried out in the same way as in Theorem
\ref{coroll}, with some additional care due to the fact that, if $f$
is the multifunction representing $S'$, in general $f(z)$ is no
longer a $C^1$ curve but only a $C^0$ one.
\enddemo
It is clear that Theorem \ref{main4} cannot be treated by applying
the methods of analytic multifunctions. As already observed, its
proof requires an accurate analysis of a single leaf.

\subsection{Preliminary results}\label{prel} First of all, we show
the following \bl \label{lochom}Let $\S$ be any complex leaf of the
foliation of $S$. Then the projection $\pi|_\S$ is a local
homeomorphism. \el \demo Let $p\in \S$, $p=(p_z,p_w)$; it suffices
to show that the differential of $\pi|\S$ is surjective in $p$. In
the opposite case, there would exist neighborhoods $U$ of $p$ in
$\C^2$ and $V$ of $p_w$ in $\C_w$, and a holomorphic function $f:V
\to \C$ such that, denoting by $\S'$ the connected component of
$\S\cap U$ which contains $p$, we would have
$$\S' = \{z=f(w)\}$$
and $\frac{\partial f}{\partial w} (p_w)=0$. In other words,
$\partial/\partial w\in T_p^\C(\S)$ and thus $\partial/\partial v
\in T_p(\S).$ This would imply
$$\frac{\partial}{\partial v}\in T_p(S),$$ which contradicts the fact that $\rho$ is a smooth
function on $\C_z\times \R_u$. \enddemo

Lemma \ref{lochom} shows that a complex leaf $\S$ of the foliation
is locally a graph over $\C_z$, but, since we do not know whether or
not $\pi:\S\to \C_z$ is actually a covering, we cannot conclude
immediately that $\pi|_\S^{-1}$ is single-valued. However, if this
is the case, it is easy to deduce that the thesis of Theorem
\ref{main7} holds true for $\S$, provided that the projection
$\pi|_\S$ is onto:

\bl \label{surj} Let $\S$ be a complex leaf of $S$, and suppose that
\begin{itemize}
\item $\pi(\S)=\C_z$;
\item for every $z_0\in \C_z$, $\pi^{-1}(z_0)\cap \S$ is a single
point.
\end{itemize}
Then there exists $c\in \C$ such that
$$\S = \{w=c\}.$$\el \demo Indeed, in this case the leaf $\S$ is biholomorphic to
$\C$ as $\pi|_\S$ is one to one; then, denoting by $v$ the
projection on the $v$-coordinate, $v\circ(\pi|_\S)^{-1}$ is a
harmonic, bounded function on $\C_z$, which is constant by
Liouville's Theorem. Therefore $v|_\S$ is also constant and so is
$u|_\S$, which is conjugate to $v$ in $\S$.
\enddemo

\br One may ask whether the latter hypothesis in Lemma \ref{surj}
can be replaced by
\begin{itemize}
\item $\pi|_\S$ is a local homeomorphism.
\end{itemize}
This is not the case, as it is shown by the following example.\er
\bex Consider the subset
$$L=\{(x,y)\in\R^2: |y|<1 \}.$$
It is simple to show that there exists a map $\phi:L\to\C$ such that
\begin{itemize}
\item $\phi$ is onto;
\item the differential of $\phi$ is always invertible;
\item $\phi$ extends as a continuous function $\oli L\to \C$;
\item $\phi^{-1}(z)$ consists of finitely many point for every
$z\in\C$.
\end{itemize}
In order to be convinced of this fact, one may proceed as follows. A
map of this kind can be identified with the smooth motion of an open
segment on $\C$ along a curve parametrized by $\R_x$. For instance,
we can first cover a ball $B\subset \C$ by this motion, and then let
the segment proceed along a suitably chosen spiral to fill in the
whole $\C$.

 Then, denoting by $J$ the standard complex structure on $\C$, we
can endow $L$ with the complex structure $\phi^\star(J)$, thus
obtaining a simply connected open Riemann surface $L_\C$ for which
$\phi:L_\C\to\C$ is tautologically holomorphic. By Riemann's
uniformization Theorem, we find a biholomorphism $\psi:L_\C\to X$
where $X$ is either $\C$ or the unit disc $D\subset\C$. We claim
that the first case cannot occur. Indeed, consider the set
$$A_\e= L_\C \cap D((0,1),\e)$$
where $D((0,1),\e)$ is a disc with center $(0,1)\in \R^2$ and radius
$\e<<1$. If $\e$ is small enough, $A_\e$ is mapped by $\phi$
biholomorphically onto an open set of $\C$. Moreover, consider $\psi
(A_\e)\subset \C\subset \C\P^1$, and observe that the boundary of
$\psi(A_\e)$ contains $\{\infty\}$ and (for small $\e$) $0\notin
\psi(A_\e)$. If $\{\mathcal U_n\}$ is a fundamental system of
neighborhoods of
$$\oli A_\e \cap \{y=1\}$$ in $\oli A_\e$, it is also clear that the
sets $\psi(\mathcal U_n)$ approach $\infty\in\C\P^1$ as $n\to
+\infty$. Let $g\in \mathcal O(\C\P^1\setminus\{0\})$ be such that
$g(\infty)=0$ and $g\not\equiv 0$. Then, in view of the previous
observations, the function $f\in \mathcal O(\phi(A_\e))$ defined by
$f=g\circ\psi\circ\phi^{-1}$ is continuous up to $\phi(\{y=1\})$,
and vanishes on this set. This is a contradiction (see also Lemma
\ref{Rado}).

It follows that $X=D$. Let $i:D\to \C^2$ be defined as
$$i(z)=(\phi\circ\psi^{-1}(z),z);$$
then $i$ is a holomorphic embedding of $D$ in $\C^2$ and we set
$\S=i(D)$. Observe that $\pi:\S\to\C$ is onto and a local
homeomorphism; moreover, by construction $|v|<1$ on $\S$. It is
clear that,for topological reasons, $\S$ cannot be the leaf of any
foliation of a graph in $\C^2$. \eex

In order to prove Theorem \ref{main7} our aim is to apply Lemma
\ref{surj} and so, from now on, we shall focus on a single complex
leaf $\S$ of the foliation of $S$ and we will prove that its
projection over $\C_z$ is a biholomorphism.

We set
$$\pi(\S)=\O\subset\C_z;$$
$\O$ is an open subset of $\C_z$.

We suppose, by contradiction, that $\O\subsetneq \C_z$, and let
$z_0\in b\O$. The following result shows that $z_0$ must actually
belong to $\O$ at least in some special case, thus proving that
$b(\O)=\emptyset$ in such a situation.

\bl \label{cluster}Let $z_0\in \C_z$ and suppose that there exist
$p_0$ such that $\pi(p_0)=z_0$ and $p_0$ is a cluster point for
$\S$. Then $z_0\in \O$. \el \br Since we do not know, at this stage,
whether or not $\S$ is a closed submanifold, it is a priori possible
that $p_0\notin\S$. Nevertheless, $\pi^{-1}(z_0)\cap\S\neq
\emptyset$.\er

\noindent \textbf{Proof of Lemma \ref{cluster}}. Let $V$ be a
neighborhood of $p_0$ on which the foliation of $S\cap V$ is
trivial.  Then, either $\S\cap V$ has finitely many connected
components - in this case one of them must contain $p_0$ - or the
connected components of $\S\cap V$ accumulate to the leaf $\S'$ of
$S\cap V$ containing $p_0$. Then $\S'$ must be a complex leaf, too.
Thus, from Lemma \ref{lochom} it follows that, if $V'\Subset V$
($p_0\in V'$) is small enough, all the leaves of $S\cap V'$
intersect (possibly in $V$) $\pi^{-1}(z_0)$. By hypothesis
$$V'\cap \S\neq\emptyset,$$
thus $\S$ contains a leaf of $S\cap V'$, therefore
$$\pi^{-1}(z_0)\cap \S\neq \emptyset.$$\enddemo

In section \ref{Shche}, applying the results of \cite{Shc}, we will
prove that $\pi|_\S^{-1}$ is single-valued. Then, given $z\in \O$,
we will denote by $w(z)$ (respectively $u(z)$,$v(z)$) the
$w$-coordinates (resp. the $u$- and $v$-coordinate) of
$\pi|_\S^{-1}(z)$. With these notations, we can state the following
straightforward corollary of Lemma \ref{cluster}:

\bc Let $z_0\in b\O$, and let $\{\mathcal U_k\}_{k\in \N}$ be a
fundamental system of neighborhoods of $z_0$ in $\C_z$. Then, for
any $M>0$ there exists $K\in\N$ such that $|w(z)|>M$ for all $z\in
\O\cap\mathcal U_k$ with $k\geq K$.\ec \demo Otherwise, there would
exist $M>0$ and a sequence $\{z_n\}_{n\in \N}$ such that
\begin{itemize}
\item $z_n\in \O$ for every $n\in \N$;
\item $z_n\to z_0$;
\item for every $n\in \N$ there exists $p_n\in \S$ such that $\pi(p_n)=z_n$ and $|w(p_n)|\leq M$.
\end{itemize}
Then $\{p_n\}_{n\in \N}$ would admit an accumulation point $p_0$ in
$\C^2$ such that $\pi(p_0)=z_0$. By Lemma \ref{cluster} this would
imply $z_0\in \O$, a contradiction. \enddemo Since, by the main
hypothesis, $v(z)$ is bounded on $\O$, it follows immediately

\bc \label{unbound}Let $z_0\in b\O$, and let $\{\mathcal U_k\}_{k\in
\N}$ be a fundamental system of neighborhoods of $z_0$ in $\C_z$.
Then for any $M>0$ there exists $K_0\in\N$ such that $|u(z)|>M$ for
all $z\in \O\cap\mathcal U_k$ with $k\geq K_0$.\ec

\br In any case, since $\O \cap \mathcal U_k$ needs not be connected
even for large $k$, it is possible that $u$ assumes both signs in
every neighborhood of $z_0$. Later on we are going to prove that it
is not the case.\er

\subsection{Unbounded harmonic functions on the disc}
Our strategy is now to reduce our situation to a problem on the
disc. For this, we will need some results about holomorphic
functions on $D=\{z\in \C: |z|<1\}$.

The following one shows that a conjugate to a bounded harmonic
function on $D$, although not necessarily bounded, cannot go to
infinity on too \lq\lq large\rq\rq\ a subset of the boundary.

\bl \label{conj}Let $f\in \mathcal O(D)$, $u= Re f$ and $v=Im f$.
Suppose that there exists a non-constant arc $\g$ in $bD$ such that
for every $M>0$ there exist a neighborhood $U_M$ of $\g$ in $\C$
such that
$$u(z)>M \ \forall z\in U_M\cap D.$$
Then $v$ is not bounded on $D$.\el \demo Taking polar coordinates
$(r,\t)$, we may assume that
$$\g=\{r=1,-\e \leq \t \leq \e\}$$
for some $\e>0$. Consider the family of arcs
$$\g_t = \{r=t, -\e \leq \t \leq \e\}$$
and set (for $t<1$)
$$I_t=\int_{\g_t}u(\t)d\t$$
the hypothesis of the Lemma implies that $I_t\to +\infty$ as $t\to
1$. Take $\oli t>0$; then, if we define (for $t>\oli t$) $J_t = I_t
- I_{\oli t}$, we have
$$\int_{\oli t\leq r \leq t, -\e \leq \t \leq \e} \frac{\partial u}{\partial r} dr d\t = \int_{-\e \leq \t \leq \e} \left ( \int_{\oli t\leq r \leq t} \frac{\partial u}{\partial r} dr \right ) d\t = $$
$$ = \int_{-\e \leq \t \leq \e} (u(t,\t) - u(\oli t,\t)) d\t = I_t - I_{\oli t} = J_t$$
and $J_t\to +\infty$ as $t\to 1$. Thus, by the integral mean value
Theorem, for all $N>0$ there exist $t'>\oli t$ such that
$$ \int_{-\e \leq \t \leq \e} \frac{\partial u}{\partial r}(t',\t) d\t > N.$$
On the other hand, in view of Cauchy-Riemann equations we have $u_r
= -v_\t$ and consequently
$$\int_{-\e \leq \t \leq \e} \frac{\partial u}{\partial r}(t',\t) d\t =
\int_{-\e \leq \t \leq \e} - \frac{\partial v}{\partial \t}(t',\t)
d\t = v(t',-\e) - v(t',\e),$$ whence $v$ is unbounded. \enddemo

\bex The previous result fails to be true if we drop any hypothesis
which guarantees that $u$ goes to infinity on a sufficiently large
subset of the boundary of $D$. A simple example is the following.
Consider the set
$$L=\{z\in\C:|y|<1\}$$
and choose a biholomorphism $\phi:L\to D$. Then
$$u=x\circ\phi^{-1}:D\to\R$$
is a harmonic function, conjugate to
$$v=y\circ\phi^{-1}:D\to\R;$$
but $u$ is not bounded while $|v|<1$ on $D$. In this case, the upper
and lower level sets of $u$ approach (two) isolated points on $bD$,
rather than segments of positive measure.\eex

\subsection{Analysis of $\O$} \label{Shche}

Our purpose now is to show that $\O$ is simply connected, which will
allow us to apply Riemann's mapping Theorem and then Lemma
\ref{conj}. In order to achieve this we apply some of the results of
\cite{Shc}, in particular the in-depth analysis which is carried out
therein on the leaves of the foliation of the Levi-flat solution for
graphs. First of all, we prove that $\pi|_\S^{-1}$ is actually
single-valued over $\O$.

\bl \label{single} Let $\O$ and $\S$ be as above. Then
$\pi|_\S^{-1}(z)$ consists of a point for every $z\in \O$. \el \demo
Suppose that, for some $z\in \O$, there exist $p,q\in \S$ ($p\neq
q$) such that $\pi(p)=\pi(q)=z$. Since, by definition, $\S$ is
connected, there exists an arc $\widetilde \g$ which joins $p$ and
$q$. Let $\g=\pi\circ\widetilde\g$ be the corresponding loop in
$\O$. Let $B$ be a ball in $\C_z\times \R_u$, centered at $z$, with
a large enough radius such that $\gamma\subset B$ and
$\tau\circ\widetilde\gamma\subset B$. Then
$$S\cap\tau^{-1}(B)=\G(\rho|_B)\subset\C^2$$
is the Levi flat surface which has the graph
$$S\cap\tau^{-1}(bB)=\G(\rho|_{bB})$$
as boundary. By the results in \cite{Shc}, we conclude that each
leaf of the foliation is properly embedded in $S\cap\tau^{-1}(B)$
(observe that, under the hypothesis of Theorem \ref{main4}, this
fact is granted by our assumption) and, therefore, that $\tau(\S)$
is properly embedded in $B$. By the choice of $B$, $\tau(p)$ and
$\tau(q)$ belong to the same connected component of $\tau (\S)\cap
B$; let $\S'$ be this component. Lemma \ref{lochom} implies that
$\S'$ is locally a graph over $\C_z$; since $B$ is convex, by virtue
of Lemma 3.2 in \cite{Shc} we deduce that $\S'$ is globally a graph
over some subdomain of $\O$. Since $\tau(p)$ and $\tau(q)$ have the
same projection over $\O$, it follows $\tau(p)=\tau(q)$ and
consequently $p=q$, a contradiction.
\enddemo

By the previous Lemma $\S$ is represented by the graph of a
holomorphic function over $\O$. Let us denote by $u$ (respectively
$v$) the real (respectively the imaginary) part of this function.
The following Lemma is an immediate consequence of the results in
\cite{Shc}:

\bl $\O$ is simply connected. \el \demo Observe that, if $\O$ is not
simply connected, then $D\cap \O$ is not simply connected for some
open disc $D\subset \C_z$. Arguing as in the previous Lemma, we
prove that $\tau(\S)$ is properly embedded on some subdomain
$$D\times(-R,R)\subset\C_z\times\R_u,\ R>>0$$
(again, under the hypothesis of Theorem \ref{main4} this is a direct
consequence of the assumption). But $\tau(\S)$ is the graph of $u$
over $D\cap \O$; since $v$ is a single-valued harmonic conjugate of
$u$, we can apply Lemma 3.3 of \cite{Shc} to obtain that $D\cap \O$
is in fact simply connected.
\enddemo

Because of the previous Lemma, we can consider a biholomorphic map
$\mathfrak R:\O\to\D$, where $D\subset\C$ is the unit disc. Our aim
is to apply Lemma \ref{conj} and in order to do so we must examine
the behavior of $\mathfrak R$ near the boundary of $\O$.

\subsection{Proof of Theorem \ref{main4}}

\bl \label{connect}Let $C$ be a connected component of the boundary
of $\O$. Then there exist a neighborhood $\mathcal U$ of $C$ in
$\oli\O$ such that either $u>0$ on $\mathcal U$ or $u<0$ on
$\mathcal U$.\el \demo Let $K$ be a compact connected subset $C$; it
is enough to prove that the thesis holds for any such $K$. Observe
that, since $\O$ is connected, $C\setminus K$ has at most two
connected components. By Corollary \ref{unbound}, for any $z\in K$
there exist a disc $D(z,\e)$ such that $|u|>0$ on $D(z,\e)\cap \O$;
thus $K$ can be covered by a finite set $\{D_1,\ldots,D_k\}$ of such
discs. If $\d$ is small enough, then
$$\mathcal U'=\{z\in \C: d(z,K)<\d\}\subset D_1\cup\ldots\cup D_k.$$
The thesis then follows from the fact that there is a connected
component of $\mathcal U'\cap \O$ whose boundary contains $K$.
Suppose that this is not the case, and choose a connected component
$\mathcal V$ of $\mathcal U'\cap \O$ such that $E= b\mathcal V\cap
K\neq \emptyset$. Observe that $b\mathcal V=E\cup F\cup G$, where
$$F=b\mathcal V \cap \{z\in \C: d(z,K)=\d\}\ {\rm and} \ G=b\mathcal
V\cap C\setminus K;$$ obviously $E\cap F=\emptyset$ and thus $G$ has
at least two connected component. Moreover, $E$ is connected since
otherwise $C\setminus K$ would have more than two connected
components. But if $E\subsetneq K$ is connected then it can touch at
most one connected component of $C\setminus K$ and thus of $G$; it
follows $E=K$. \enddemo

\bc Let $C$ be a connected component of $b\O$. Then there is a
fundamental system $\{\mathcal V_n\}_{n\in\N}$ of neighborhoods of
$C$ in $\oli \O$ such that either $$\inf_{\mathcal V_n}u\to
+\infty$$ or $$\sup_{\mathcal V_n}u\to -\infty$$ as $n\to \infty$.
\ec \demo This is a consequence of Corollary \ref{unbound} and Lemma
\ref{connect}. \enddemo

\br \label{stupid} In the previous statement, we can assume that
$\mathcal V_n$ and $\O\setminus\mathcal V_n$ are connected for every
$n$. We can also assume that the sequence $\{\mathcal
V_n\}_{n\in\N}$ is decreasing, with $\oli {\mathcal V}_{n+1}\subset
\mathcal V_n$ (where the closure is taken in $\O$).\er Now we fix
our attention to the sequence $\{\mathcal W_n= \mathfrak R(\mathcal
V_n)\}_{n\in\N}$ of domains of $D$. For each $n\in N$, we define
$$\L_n = \oli{\mathcal W}_n \cap bD$$
where the closure of $\mathcal W_n$ is taken in $\C$.

\bl $\{\L_n\}_{n\in\N}$ is a non-increasing sequence of closed,
connected subsets of $bD$; moreover, $\L_n\neq\emptyset$ for every
$n\in \N$. \el \demo We prove various points separately:
\begin{itemize}
\item each $\L_n$ is closed by definition, and the sequence is
non-increasing by Remark \ref{stupid};
\item $\L_n\neq \emptyset$ for, in the opposite case, we would have
$\mathcal W_n\Subset D$, which implies that $\mathfrak R$ is not
proper. This is a contradiction because $\mathfrak R$ is a
biholomorphism;
\item $\L_n$ is connected because otherwise we would have that
$$D\setminus \mathcal W_n=\mathfrak R(\O \setminus \mathcal V_n)$$ is
not connected, which would contradict the choice of $\mathcal V_n$
(see Remark \ref{stupid}). Indeed, suppose that $\L'$ and $\L''$ are
two disjoint connected components of $\L_n$, and choose a simple arc
$\g\subset \mathcal W_n$ joining a point of $\L'$ and a point of
$\L''$. Then $D\setminus \g$ has two connected components $D_1$ and
$D_2$, and we must have $D_j\cap (D\setminus \mathcal W_n)\neq
\emptyset$ for $j=1,2$ since $\L'$ and $\L''$ are disconnected. This
implies that $D\setminus \mathcal W_n$ is disconnected.
\end{itemize} \enddemo
From the previous Lemma it follows that, if we set
$$\L=\bigcap_{n\in\N}\L_n$$
then $\L$ is a closed, non-empty interval of $bD$. A priori, $\L$
could be reduced to a single point of $bD$. In order to apply Lemma
\ref{conj}, we must prove that this is not the case. \bl
\label{Rado} The interval $\L\subset bD$ is not reduced to a
point.\el \demo We argue by contradiction and assume that
$\L=\{z_0\}$ with $z_0\in bD$. Observe that, for every $\e>0$, there
exists $N\in\N$ such that
$$\mathcal W_n\subset D(z_0,\e)\cap D\ \forall n\geq N,$$
where $D(z_0,\e)$ is the disc centered at $z_0$ with radius $\e$.
Indeed, in the opposite case there would exist (cfr. Remark
\ref{stupid}) $p\in D$ such that $p\in \mathcal W_n$ for all $n\in
\N$, and this is not possible since $\mathcal V_n$ is a fundamental
system of neighborhoods of $C$. Now consider, on $D$, the
holomorphic function $f(z)=z-z_0$, and let $g\in \mathcal O(\O)$ be
defined by $g = f\circ \mathfrak R$. Then, by the choice of $f$ and
the previous observation, $g$ extends to $\O\cup C$ as a continuous
function putting $g\equiv 0$ on $C$. Choose a point $w\in C$ and
consider a disc $D'=D(w,\e)$ such that $D'\setminus C$ is
disconnected. Define a function $\widetilde g:D'\to \C$ by
$$\widetilde g(z)=\left\{
                 \begin{array}{ll}
                   g(z), & \hbox{$z\in \oli\O\cap D'$;} \\
                   0, & \hbox{$z\in D'\setminus \oli\O$.}
                 \end{array}
               \right.
$$
Then $\widetilde g$ is continuous. Moreover, by definition
$\widetilde g$ is holomorphic outside the set $\{\widetilde g=0\}$;
therefore, by Rado's Theorem, $\widetilde g\in \mathcal O(D')$.
Since $\stackrel{\circ}{\{\widetilde g=0\}}\neq \emptyset$, we have
$\widetilde g\equiv 0$ on $D'$ and consequently $g\equiv 0$ on $\O$,
which is a contradiction. \enddemo Now we are in position to prove
Theorem \ref{main4}: Lemma \ref{Rado} allows us to apply Lemma
\ref{conj} and deduce that $u$ cannot be unbounded on $\O$. By
Corollary \ref{unbound} we have that $\O=\C_z$ and thus $\pi$ is
onto. Lemma \ref{single} implies that $\pi$ is one to one, therefore
we can apply Lemma \ref{surj} and conclude that $\S=\{w=c\}$ for
some $c\in\C$, whence the thesis of Theorem \ref{main7}. \enddemo

\subsection{The result in $\C^n$} The statement of Theorem
\ref{main7} can be generalized to the case when $S$ is a Levi-flat
hypersurface of $\C^n$. Consider holomorphic coordinates
$(z_1,\ldots,z_{n-1},w)=(z,w)$, $z_j = x_j + i y_j$, $w=u+iv$, and
let $\rho:\C^{n-1}\times \R \to \R$ be a smooth function such that
$S=\{v=\rho(z,u)\}$ is a Levi-flat graph. Then we can restate almost
verbatim Theorem \ref{main7}: \bt \label{main2}$S$ is foliated by
complex hyperplanes, i.e. $\rho$ has the form $$\rho (z,u) =
\varrho(u)$$ for some smooth function $\varrho:\R\to\R$. \et \demo
This is an easy consequence of Theorem \ref{main7}. Indeed, let
$p_1=(z_1,u)$ and $p_2=(z_2,u)$ be two points in $\C^{n-1}_z\times
\R_u$ with the same $u$-coordinate, and consider the complex line
$L\subset \C^{n-1}_z$ such that $z_1,z_2\in L$. Then the restriction
of $\rho$ to $L\times\R_u$ has a Levi-flat graph
$$S_L=S\cap (L\times \C_w) \subset L\times\C_w\cong \C^2.$$
Theorem \ref{main7} applies to $S_L$, showing that $\rho|_{L\times
\R_u}$ is a function of $u$ and thus that $\rho(p_1)=\rho(p_2)$.
This proves the thesis.
\enddemo

\subsection{Generalization to a continuous graph}

The arguments of the previous sections work in the case that $\rho$
is at least of class $C^2$. However, it is possible to generalize
the result to the case of a continuous graph. In order to achieve
this, the Main Theorem of Shcherbina's paper \cite{Shc} (which gives
also a description of the leaves of the foliation of the polynomial
hull of a graph in $\C^2$) can be applied, rather than Lemmas 3.2
and 3.3. We say that a continuous hypersurface $S\subset\C^n$ (i.e.
a subset which is locally a graph of a continuous function over an
open subset of a real hyperplane of $\C^n$) is \emph{Levi flat} if
it (locally) separates $\C^n$ in two pseudoconvex domains. Note
that, in the case $n=2$, $S$ is locally the union of a disjoint
family of complex discs (see again \cite{Shc}, Corollary 1.1). So,
let $\rho:\C^{n-1}\times \R \to \R$ be a continuous function such
that $S=\{v=\rho(z,u)\}$ is a Levi flat graph. Then, as before, we
have \bt \label{main3} $S$ is foliated by complex hyperplanes, i.e.
$\rho$ has the form $$\rho (z,u) = \varrho(u)$$ for some continuous
function $\varrho:\R\to\R$. \et Once again it is sufficient to show
that the statement is true for $n=2$. In this case we do not know a
priori whether $S$ has a foliated atlas; nevertheless, since each
$p\in S$ is contained in a germ of holomorphic curve $\S_p\subset S$
(and this germ is unique in view of Lemma 4.1 of \cite{Shc}) we can
still consider the maximal connected surface $\S$ that passes
through $p$. Our aim is to carry out an analysis of $\S$ similar to
that delivered in the previous sections for the $C^2$ case. First of
all, we want to generalize Lemma \ref{lochom}: \bl
\label{lochom2}Let $\S$ be any leaf of the foliation of $S$. Then
the projection $\pi|_\S$ is a local homeomorphism. \el \demo In this
case the fact that $\partial/\partial v\in T(\S)$ does not give a
contradiction, since $S$ is only a continuous graph. Instead, we
rely on the Main Theorem of \cite{Shc} Let $p\in \S$, $p=(p_1,p_2)$,
$B$ a ball in $\C_z\times\R_u$ containing the point $(p_1,\Re p_2)$
and consider $\rho|_{bB}$. Then Shcherbina's Theorem applies to
$\gamma = \Gamma(\rho|_{bB})$, hence by the point (ii) of that
statement it follows that the disc through $p$ is a graph over a
domain of $\C_z$.\enddemo

As before, we define $\O=\pi(\S)$ and we prove that $\pi|_\S^{-1}$
is single-valued and that $\O$ is simply connected.

\bl \label{single2} Let $\O$ and $\S$ be as above. Then
$\pi|_\S^{-1}(z)$ consists of a point for every $z\in \O$. \el \demo
We follow the proof of Lemma \ref{single}. Suppose that, for some
$z\in \O$, there exist $p,q\in \S$ ($p\neq q$) such that
$\pi(p)=\pi(q)=z$. We choose an arc $\widetilde \g$ joining $p$ and
$q$, with $\g=\pi\circ\widetilde\g$ the corresponding loop in $\O$.
Let $B$ be a ball in $\C_z\times \R_u$, centered at $z$, with a
large enough radius such that $\gamma\subset B$ and
$\tau\circ\widetilde\gamma\subset B$. Then
$$S\cap\tau^{-1}(B)=\G(\rho|_B)\subset\C^2$$
is the Levi-flat surface which has the graph
$$S\cap\tau^{-1}(bB)=\G(\rho|_{bB})$$
as boundary. Since, by Shcherbina's Main Theorem, $S\cap
\tau^{-1}(B)$ is the disjoint union of discs which are graphs on
$\C_z$, we must have $p=q$, which is a contradiction.
\enddemo

\bl $\O$ is simply connected. \el \demo As in the previous case, we
assume by contradiction that $D\cap \O$ is not simply connected for
some open disc $D\subset \C_z$. We choose a ball $B\subset
\C_z\times\R_u$ such that $B\cap \C_z=D$. Then, by point (ii) of
Shcherbina's Main Theorem, the leaves of $S\cap \tau^{-1}(B)$ are
graphs over simply connected domains of $\C_z$. It follows that
$$D\cap \O=\pi (\S\cap \tau^{-1}(B))$$ must be simply connected.
\enddemo

Now we prove the analogous of Lemma \ref{cluster}:

\bl \label{cluster2}Let $z_0\in \C_z$ and suppose that there exist
$p_0$ such that $\pi(p_0)=z_0$ and $p_0$ is a cluster point for
$\S$. Then $z_0\in \O$. \el \demo In this case we can actually prove
that $p_0\in \S$, i.e. $\S$ is a closed surface. Ideed, consider a
ball $B\subset \C_z\times \R_u$ which is centered at $\tau(p_0)$.
Then $S_B=S\cap\tau^{-1}(B)=\G(\rho|_B)$ is a union of disjoint
complex discs which are graphs over domains of $\C_z$. Since $\S$ is
a graph over $\C_z$ and contains points of $S_B$, it must contain
exactly one of those discs, which has to be the one passing through
$p_0$, $p_0$ being a cluster point. Then $p_0\in\S$.
\enddemo
Keeping the notation adopted in section \ref{prel}, we then have,
with the same proof as \ref{unbound}, the following \bc
\label{unbound2}Let $z_0\in b\O$, and let $\{\mathcal U_k\}_{k\in
\N}$ be a fundamental system of neighborhoods of $z_0$ in $\C_z$.
Then for any $M>0$ there exists $K_0\in\N$ such that $|u(z)|>M$ for
all $z\in \O\cap\mathcal U_k$ with $k\geq K_0$.\ec

The rest of the proof of Theorem \ref{main3} goes exactly as in the
previous case.

\section{Foliations of $D\times \C$} \label{folcil}

Let $D\subset \C$ be the unit disc. As seen in the section
\ref{uno}, the methods of analytic multifunctions allow to prove a
\lq\lq Liouville result\rq\rq\ for Levi flat manifolds contained in
$D\times \C$ by considering them \lq\lq as a whole\rq\rq. In such a
case, since the leaves of the foliation are uniquely determined, we
obtain immediately that the these leaves are complex lines.  In what
follows, we want to consider smooth foliations of $D\times \C$ by
complex curves. In this case the global object defines a constant
multifunction, but this does not imply that the foliation is
trivial. So, in order to show that - under natural topological
restriction - this is the case, we need to study the foliation
\lq\lq leaf-by-leaf\rq\rq\ as done in section \ref{folgraf}.

We observe that, even though The kind of foliations we treat is not
perhaps the most general case in which a result of Liouville's type
can be obtained, nevertheless the methods of multifunctions do not
apply directly.

Choose holomorphic coordinates $(z,w)$ in $\C^2$ such that $D$ is
the unit disc on $\C_z$: \bd \label{regolare} We say that a
$2$-dimensional real foliation of $D_z\times \C_w$ is \emph{regular}
if, for any leaf $\S$ of the foliation, the following conditions are
fulfilled:
\begin{description}
  \item[(i)] for any open subset $U\Subset \C_w$, every connected component of $S\cap (D_z\times
  U)$ is (up to a possible shrinking of $U$) a compact $2$-manifold
  with boundary which has positive distance from $bD_z\times \C_w$;
  \item[(ii)] there exists a ball $B\Subset \C_w$ such that every connected component of
  $$\S\cap (D_z\times (\C_w\setminus B))$$ is finitely branched over
  $D_z$.
\end{description} \ed
Roughly speaking, condition (i) says that the foliation has to be
defined only on the cylinder, in such a way that it is not possible
to \lq\lq extend it further\rq\rq. Condition (ii) guarantees a nice
behavior at infinity. A leaf is allowed to approach the boundary of
the cylinder, but only along paths which are not contained in any
compact set.

We say that a $2$-dimensional foliation on $D_z\times \C_w$ is
\emph{trivial} if all the leaves are of the form $\{z=c\}$.

\bt \label{cilindro} Let $\mathcal F$ be a (real) $2$-dimensional
regular foliation on $D\times \C$. Suppose that all the leaves of
$\mathcal F$ are complex manifolds. Then $\mathcal F$ is trivial.\et

In order to prove the Theorem we concentrate on a single leaf and
prove that it must be of the form $\{z=c\}$. So, let $\S$ be a leaf
of $\mathcal F$. We have the following

\bl \label{sotto2} Let $L\neq \C_z$ be a complex line, and denote by
$\pi$ the orthogonal projection over $L$. Then $\pi|_\S$ is onto.
\el \demo We shall argue as in the proof of Lemma \ref{sotto}. Since
$\pi|_\S$ is a holomorphic function we have that $\pi(\S)=\O$ is an
open subset of $L$. Let $p\in b\O$, and let $B$ be any ball of $L$
centered at $p$. Then, by the hypothesis $L\neq \C_z$, defining
$$U=\pi^{-1}(B\cap \O)\cap (D\times \C)$$ we have $U\Subset \C^2$ and thus the connected
components of $\S\cap U$ have positive distance from $bD\times \C$.
Let $\S'$ be one of these components; then $\pi(\S')=\O\cap B$.
Otherwise, we would have $\pi(S')=\O'$ where $\O'\subset \O\cap B$
is an open subset. Take $q\in b\O'\setminus b(\O\cap B)$; then,
arguing exactly as in the proof of Lemma \ref{sotto} (taking in
account the fact that, since $L\neq \C_z$, any complex line which is
orthogonal to $L$ and contains points of $D\times \C$ touches
$bD\times \C$), we have $q\in \pi(\S')$. It follows that $p\in
\oli{\pi(\S')}$. Arguing again as in \ref{sotto} we conclude that
$p\in \pi(\S')$, hence $p\in \pi(\S)$.
\enddemo

The previous Lemma suggests to consider $\S$ as an analytic
multifunction defined over $\C_w$ (or any complex line $L\neq
\C_z$). However, while $\oli \S$ is actually an analytic
multifunction, it may have non-empty interior: $\S$ may be even
dense on $D\times \C$, in which case the fact that $\oli\S$ is a
constant multifunction gives no information.

So, we proceed as follows. Let $\widetilde \S$ be the universal
covering of $\S$. By Riemann's uniformization Theorem, $\widetilde
\S$ is either $\C$ or the unit disc $D$. In the first case the
thesis would follow immediately, since the lift of $z|_\S$ to
$\widetilde S$ would be a bounded holomorphic function on
$\widetilde \S \cong \C$, hence constant. Then, let us suppose
$\widetilde \S \cong D$. Let $f_1,f_2:D \to \C$ be defined in such a
way that $f=(f_1,f_2):D\to \C^2$ is the covering map $D\to \S$. Then
$f_1,f_2$ are holomorphic functions, $f_1$ is bounded (by $1$) and,
by Lemma \ref{sotto2}, $f_2$ is onto.

We consider a complex parametrization of $\P(\C^2)\setminus \C_w$
(i.e. the complex linear subsets different from $\C_w$) given by
$v_\eta=(\eta,1)$, $\eta\in \C$. We project $\S$ along each of the
lines of this parametrization, i.e. we consider, for any $\zeta\in
D$,
$$\left \langle f(\zeta),v_\eta^\perp \right \rangle =
\left \langle \left (f_1(\zeta),f_2(\zeta)\right ),\left
(1,-\oli\eta \right ) \right \rangle=
 f_1(\zeta) - \eta f_2(\zeta).$$
Denote by $F$ the function $D_\zeta\times \C_\eta\to \C$ defined by
the previous expression, i.e. $F(\zeta,\eta) = f_1(\zeta) - \eta
f_2(\zeta)$. Then $F$ is a holomorphic function, whose zero locus is
a $1$-dimensional analytic subset $\mathcal Z$ of $D_\zeta\times
\C_\eta$. $\oli{\mathcal Z}$ is an analytic multifunction. Since
$$\oli{\mathcal Z}\cap (D\times \C_z)=\mathcal Z$$ the property that $\oli{\mathcal Z}$ is
a constant multifunction would imply that $\mathcal Z$ is constant,
i.e. union of complex lines (which is not true for $\S$).
Unfortunately, it can occur that $\oli{\mathcal Z}$ contains
$bD_\zeta\times \C_\eta$. In such a case we would have
$\widehat{\oli{\mathcal Z}}=D_\zeta\times \C_\eta$ (that is,
$\widehat{\oli{\mathcal Z}}$ is a constant multifunction as
expected), but we could not conclude anything about $\mathcal Z$.
The difficulty behind this obstacle is that, in this case, it is not
sufficient to \lq\lq test\rq\rq \ the behavior of $\mathcal Z$ by
plurisubharmonic functions defined in a neighborhood of
$D_\zeta\times \C_\eta$, since they can \lq\lq detect what
happens\rq\rq \ only in their maximum sets, i.e., in our case,
$bD_\zeta\times \C_\eta$. Then, we need to analyze $\mathcal Z$ in
more detail.

Let $\mathcal A\subset D_\zeta$ be the subset
$$\mathcal A = \{\zeta\in D: f_1(\zeta)=f_2(\zeta)=0\};$$
then $\mathcal A$ is a discrete subset of $D_\zeta$. Observe that
$\mathcal Z$ can be expressed as
$$\mathcal Z = \{(\zeta,\eta)\in D\times \C: \zeta \in \mathcal A\}
\cup \{(\zeta,\eta)\in D\times \C: \eta=f_1(\zeta)/f_2(\zeta)\}$$
i.e. $\mathcal Z$ is the union of a discrete family of complex lines
and the graph of a meromorphic function over $D_\zeta$. Denote by
$\mathcal Z'$ this graph. Clearly $\mathcal Z'$ is a $1$-dimensional
complex submanifold of $$D_\zeta \times \C_\eta \setminus
\{(\zeta,\eta): f_2(\zeta)=0\},$$ and possibly extends as a
submanifold of a bigger domain (it extends through the complex lines
$\{\zeta=c\}$ for which $c\in \mathcal A$ and the order of zero of
$f_1$ at $c$ is bigger than or equal to the one of $f_2$). We then
take as $\mathcal Z'$ the maximal possible extension on a subdomain
of $D_\zeta\times \C_\eta$.

\bl \label{cambio} Up to a change of coordinates on $\C^2 =
\C_z\times \C_w$, for the resulting $F$ and $\mathcal Z'$ the
following is true:
\begin{itemize}
  \item $\frac{\partial}{\partial \zeta} F(\zeta,0)= \frac{\partial}{\partial \zeta} f_1(\zeta)\neq 0$;
  \item there exists a small ball $B_\e\subset \C_\eta$, centered at $0$,
  such that every connected component of $\mathcal Z'\cap
  (D_\zeta\times B_\e)$ intersects $\eta=0$;
  \item every one of such connected components is a compact
  manifold with boundary, which is a finite branched covering of $B_\e$.
\end{itemize}
 \el \demo
 We perform a coordinate change on $D_z$ such that the line $\{z=0\}$
 intersects $\S$ transversally at some point (this coordinate change
 exists, since otherwise we would have that $\S=\{z=c\}$ for some
 $c\in D_z$); in this way the first condition is assured.
 Then, we choose $\e$ small
 enough in such a way that, for any $\eta \in B_\e$, we have that
 $$L_\eta\cap (bD_z\times \C_w)\subset \C^2\setminus (B\times D_z),$$
 where $L_\eta$ is the complex line of $\C_z\times \C_w$ parametrized by
 $\eta$ and $B$ is the one in Definition \ref{regolare}.
 The choice of $\e$ implies that the last two assertions are satisfied. In
 fact, Definition \ref{regolare} implies that, for each connected
 component $\mathcal W$ of $\mathcal Z'\cap
  (D_\zeta\times B_\e)$, the fibres of the projection $\pi:\mathcal W\to
  B_\e$ are finite. Moreover, from our construction it follows that
  $\pi(\mathcal W)$ includes $B_\e \setminus \{0\}$ (this is a
  consequence of the following two facts: the intersection of $\S$ with
  $L_\eta$ is stable, and $L_\eta \cap (D_z\times
  \C_ w)$ is compact for $\eta\neq 0$). Since $\mathcal W$ is a
  finitely ramified covering of $B_\e\setminus \{0\}$, it extends to
  $D_\zeta \times B_\e$. Indeed, to see this it is sufficient to observe that
  the symmetric functions on the $\zeta$-coordinates of the fibres
  of $\pi:\mathcal W\to B_\e\setminus \{0\}$ are bounded, holomorphic
  function on $B_\e\setminus \{0\}$, thus they extend to $B_\e$.
  Since $\mathcal Z'$ is closed, $\oli{\mathcal W}\subset \mathcal
  Z'$ and therefore the connected components of $\mathcal Z'\cap
  (D_\zeta\times B_\e)$ satisfy the last two conditions of the Lemma.
\enddemo

Now we define
$$G(\zeta,\eta) = F(\zeta,\eta) - F(0,\eta) = f_1(\zeta) - \eta(f_2(\zeta) - f_2(0))$$
and denote again by $\mathcal Z'$ the graph part of the zero locus
of $G$. Observe that the arguments of Lemma \ref{cambio} work for
$G$ as well, and that $G$ satisfies
\begin{itemize}
\item $G(0,\eta)\equiv 0$
\item $\frac{\partial}{\partial \zeta}G(\zeta,\eta) \neq 0$ for
$\eta$ in a neighborhood of $0$;
\item for $\eta\neq 0$, $G(\cdot,\eta):D\to \C$ is onto;
\item for $\eta=0$, $G(\cdot,0)=f_1$ is bounded (by $1$).
\end{itemize}
We will show that these properties combined with the conclusions of
Lemma \ref{cambio} give a contradiction.

In order to do this, we need two intermediate Lemmas on the
holomorphic functions on $D$ which are consequences of the classical
Schwarz Lemma.

\subsection{Lemmas on holomorphic functions on $D$}
Let $f\in \mathcal O(D)$, $f\not\equiv 0$. Then the zero locus of
$f$ is a countable, discrete subset of $D$. We set
$$f^{-1}(0)=(a_0,a_1,\ldots,a_n,\ldots)$$
where\begin{itemize}
       \item $f(a_i)=0$ for all $i\in \mathbb N$;
       \item all the zeroes are listed except $0$ if $f(0)=0$;
       \item the $a_i$'s are listed by non-decreasing modulus;
       \item if $k$ is the multiplicity of $a_i$, then $a_i$ is
       listed $k$ times.
     \end{itemize}
Obviously, the infinite product
$$\prod_{i=1}^\infty |a_i|$$
converges to some non-negative real number. We denote it by $\Pi_0
f$. Our aim is to study the behavior of $\Pi_0 f$ in some cases,
first of all, when $f:D\to \C$ is onto.

We say that $f$ is \emph{semi-proper} if, for any compact subset
$K\subset \C$, every connected component of $f^{-1}(K)$ is a compact
subset of $D$. Every proper map is semi-proper; however, in our
setting the first property is not relevant, since there is no proper
holomorphic map from $D$ onto $\C$. \bl \label{suriet}Let $f:D\to
\C$ be a holomorphic, surjective semi-proper map. Suppose that
$f(0)=0$ and $f'(0)\neq 0$. Then $\Pi_0 f=0$.\el \demo Choose a ball
$B=B(0,R)\subset \C$, centered at $0$ and with radius $R>>0$. Let
$C$ be the connected component of $f^{-1}(\oli B)$ containing $0\in
D$; by hypothesis $C$ is a compact subset of $D$. It is easily seen
that $f:C\to \oli B$ is a finite ramified covering of $\oli B$, thus
for any $z\in B$ we may consider the finite set $f^{-1}(z)\cap
C=\{w_1(z),\ldots,w_s(z)\}$. Define $g:B\to D$ as
$$g(z)=w_1(z)\cdot\ldots\cdot w_s(z).$$
Then $g$ is a well defined holomorphic function on $B$; moreover, by
hypothesis $g(0)=0$. In view of Schwarz's Lemma,
\begin{equation}\label{nonnec}
|g'(0)|\leq 1/R.
\end{equation} Now observe that, by hypothesis, $f$
is a local homeomorphism near $0\in D$. Therefore, in a small
neighborhood $U$ of $0$ in $\C$, $g$ can be written as $g=g_1\cdot
g_2$, where $g_1$ is a local inverse of $f$ such that $g_1(0)=0$.
Taking the first derivative, we have
$$g'(0)= g'_1(0)g_2(0) + g_1(0)g'_2(0) = g'_1(0) \prod\{w: w\in f^{-1}(0)\cap C\setminus
\{0\}\}.$$ Let $k=|g_1'(0)|$; by hypothesis $k\neq 0$ and, moreover,
it clearly does not depend on the choice of $R$. By (\ref{nonnec})
we obtain
$$\left | \prod\{w: w\in f^{-1}(0)\cap C\setminus
\{0\}\}  \right | \leq \frac{1}{kR}.$$ For $R\to \infty$ we find
that the product of the modulus of a (possibly proper) subset of
$f^{-1}(0)\setminus \{0\}$ is vanishing; hence, a fortiori, $\Pi_0
f=0$. \enddemo The following result is in some sense the counterpart
of Lemma \ref{suriet}. It shows that if $f$ is bounded and
$f'(0)\neq 0$ then $\Pi_0 f\neq 0$. Indeed,

\bl \label{Schwa} Let $f\in \mathcal O(D)$ be a bounded holomorphic
function on the unit disc. Suppose that $f(0)=0$ and $\Pi_0 f = 0$.
Then $f'(0)=0$. \el \demo Choose $M>0$ such that $|f(z)|\leq M$ for
all $z\in D$. We first prove the Lemma in the case that $f\in
\mathcal O(D)\cap C^0(\oli D)$. For any $a_i\in
f^{-1}(0)\setminus\{0\}$, we choose a holomorphic automorphism
$\phi_{a_i}$ of the disc, of the form
$$\phi_{a_1}= \frac{\alpha_i z + \beta_i}{\oli\beta_i z - \oli \alpha_i}, \alpha_i,\beta_i\in \C,$$
in such a way that $\phi_{a_i}(a_i)=0$. Observe that the following
properties hold:\begin{itemize}
                  \item the vanishing order of $\phi_{a_i}$ at $a_i$
                  is $1$;
                  \item $|\phi_{a_i}(z)| = 1$ for $z\in bD$;
                  \item $|\phi_{a_i}(0)| = |a_i|$.
                \end{itemize}
Choose $\e>0$ and take $N\in \mathbb N$ big enough such that
$$|a_1|\cdot |a_2|\cdot\ldots\cdot |a_N|\leq \e.$$ Define a holomorphic function $g$ in
the following way:
$$g(z)=\frac{f(z)}{\phi_{a_1}(z)\cdot \phi_{a_2}(z)\cdot\ldots\cdot
\phi_{a_N}(z)};$$ $g$ is well defined and holomorphic in $D$ because
$\phi_{a_i}$ vanishes only at $a_i$ (of order $1$). Moreover,
$$|g(z)|= \frac{|f(z)|}{|\phi_{a_1}(z)|\cdot
|\phi_{a_2}(z)|\cdot\ldots\cdot |\phi_{a_N}(z)|} = |f(z)| \leq M \
\forall z\in bD,$$ thus $|g|\leq M$ on the whole disc $D$. Since
$g(0)=0$, by Schwarz's Lemma we obtain that $|g'(0)|\leq M$. Then,
taking the first derivative of $f$ and defining $$\Phi_N =
\phi_{a_1}\cdot \phi_{a_2}\cdot\ldots\cdot \phi_{a_N}$$ we have
$$|f'(0)|=|g'(0)\Phi_N(0)+g(0)\Phi_N'(0)|=|g'(0)|\cdot|\phi_{a_1}(0)\cdot \phi_{a_2}(0)
\cdot\ldots\cdot \phi_{a_N}(0)|=$$
$$= |g'(0)|\cdot (|a_1|\cdot |a_2|\cdot\ldots\cdot |a_N|)\leq
M\e.$$ Then, letting $\e\to 0$ we obtain the thesis when $f$ is
continuous up to the boundary. The general case is obtained by
applying the same proof to $f_\d(z)=f((1-\d)z)$ and letting $\d\to
0$.
\enddemo
\br The proof of the previous Lemma is analogous to the one of the
classical Schwarz lemma: indeed, the method is to get rid of the
zeroes of $f$ by dividing by a suitable holomorphic function $\Phi$,
and then apply the maximum principle. However, the choice of $\Phi$
needs (a bit of) care, since, for example, dividing by $1/(z-a_i)$
does not allow to obtain the right estimate. \er

\subsection{Proof of Theorem \ref{cilindro}}
Keeping the notations introduced in the previous sections, let
$$\Pi(\eta)=\Pi_0 G(\cdot,\eta).$$
As already observed, by Lemma \ref{sotto2} $G(\cdot,\eta):D\to\C$ is
onto for $\eta\neq 0$; moreover, condition (i) in definition
\ref{regolare} implies that $G(\cdot,\eta)$ is semi-proper. Hence
Lemma \ref{suriet} applies to $G(\cdot,\eta)$, showing that
$$\Pi(\eta)=0 \ \forall \eta\neq 0.$$ If we prove that $\Pi(0)=0$,
then we are in position to apply Lemma \ref{Schwa} and obtain that
$f_1'(0)=0$ (where $f_1$ is the first component of the covering map
$D\to \S$ introduced earlier), which is a contradiction. So, our
purpose is to show that $\Pi(0)=0$, and for this that $\Pi(\eta)$ is
continuous at $0\in \C_\eta$.

Observe that by Lemma \ref{cambio} we have that $\mathcal Z'\cap
(D_\zeta\times B_\e)$ is made up, for small $\e$, of countable many
connected components, each of them being a finite ramified covering
of $B_\e$. For each one of those connected components $\mathcal
Z'_i$, we define $g_i\in \mathcal O(B_\e)$ as
$$g_i(\eta)=\prod_{(\zeta,\eta)\in \mathcal Z'_i} \zeta;$$
clearly
$$\Pi(\eta)=\prod_{i=1}^\infty |g_i(\eta)|$$
for $\eta\in B_\e$. The thesis is then a consequence of the
following Lemma:

\bl \label{ultimo}Let $g_i\in \mathcal O(B_\e)$ be defined as
before. Then the product of the $g_i$'s is continuous in
$\eta=0$.\el \demo Denote by $G_k$ the product of the first $k$
functions, $$G_k(\eta)= g_1(\eta)\cdot g_2(\eta)\cdot\ldots\cdot
g_k(\eta).$$ Since $|g_i|<1$ on $B_\e$, the sequence of functions
$|G_k|$ is monotone decreasing. Moreover, the sequence $G_k$ is
uniformly bounded, hence by Montel's Theorem there is a subsequence
$G_{k_j}$ which converges uniformly to a continuous (holomorphic)
function $G$. Since we already know that $\Pi(\eta)=0$ for $\eta\neq
0$, it follows that $G\equiv 0$ on $B_\e$. Therefore, since $|G_k|$
is a decreasing sequence which admits a subsequence convergent to
zero, we conclude that $|G_k|\to 0$ (uniformly in $\eta$), i.e.
$$\Pi(0)=\prod_{i=1}^\infty |g_i(0)| = 0.$$
\enddemo

\br Part (ii) of Definition \ref{regolare} is used only in Lemma
\ref{cambio}, to assure that in a neighborhood of $\{\eta=0\}$ the
connected components of the zero locus $\mathcal Z'$ are
well-behaved, i.e.\ they are a finite branched covering of a
neighborhood of $0$ in $\C_\eta$. This allows us to prove Lemma
\ref{ultimo}. We conjecture that (ii) is superfluous and the
triviality result is still valid only under assumption (i) of
Definition \ref{regolare}; however, our method does not work in this
generality. This is due to the fact that, in general, a connected
Riemann surface which is a branched covering of a neighborhood of
$0$ in $\C_\eta\setminus \{0\}$ may be not extendable through
$\{\eta=0\}$, even if $\zeta$-coordinate of its points is bounded
(see Example \ref{sotto3} below). In such a case the proof of Lemma
\ref{ultimo} does not work. \er \bex \label{sotto3} Consider, in
$\C^2\setminus \{z=0\}$, the set defined $S$ by
$$w = \sum_{j=1}^\infty \frac{1}{2^j}\sqrt{z-\frac{1}{j}}.$$
Here, we mean the following: for each $z\neq 0$ and each $j\in \N$,
we order arbitrarily the two roots $r^j_1$, $r^j_2$ of $\sqrt{z -
1/j}$. For any function $c:\N\to \{0,1\}$, we set
$$w(c) = \sum_{j=1}^\infty \frac{1}{2^j}r^j_{c(j)};$$
this sum converges because $|r^j_{c(j)}|$ is bounded. Then $S$ is
the collection of the points $(w(c),z)$ for all $z\in
\C\setminus\{0\}$ and for all $c:\N\to\{0,1\}$. This set is not a
Riemann surface, since the fiber over $z\in \C\setminus\{0\}$ is not
countable. However, any connected component $S'$ of $S$ is a Riemann
surface, with branching points $z=1/j$, $j\in \N$, bounded in a
neighborhood of $\{z=0\}$. Clearly $S'$ does not extend to through
$\{\eta=0\}$. \eex


\begin{thebibliography}{8}
\bibitem{BI} David E. Barret and Takashi Inaba, \emph{On the topology of compact smooth three-dimensional Levi-flat
   hypersurfaces}, J. Geom. Anal. 2 (1992), no. 6, 489-497
\bibitem{Bry} Robert L. Bryant, \emph{Levi-flat minimal hypersurfaces in two-dimensional complex space
   forms}, hundred years after Sophus Lie (Kyoto/Nara, 1999) Adv. Stud. Pure
   Math., vol. 37, Math. Soc. Japan, Tokyo, 2002, pp. 1-44.
\bibitem{LT} Guido Lupacciolu and Giuseppe Tomassini, \emph{An extension Theorem for
CR-functions}, Ann. Mat. Pura Appl. (4) 137 (1984), 257-263
(Italian, with English summary)
\bibitem{Oka} Kiyosi Oka, \emph{Note sur les familles de fonctions analytiques
multiformes}, J. Sci. Hiroshima Univ. Ser. A 4 (1934) 93-98
\bibitem{Ran} Thomas Ransford, \emph{A new approach to analytic
multifunctions}, Set-Valued Anal. 7 (1999), no. 2, 159-194
\bibitem{Shc} Nikolay V. Shcherbina, \emph{On the polynomial hull of a
graph}, Indiana Univ. Math. J. 42 (1993), no. 2, 477-503
\bibitem{Slo} Zbigniew Slodkowski, \emph{Analytic set-valued functions and
spectra}, Math. Ann. 256 (1981), no. 3, 363-386
\bibitem{Sto} Gabriel Stolzenberg, \emph{Uniform approximation on smooth
curves}, Acta Math. 115 (1966), 185-198
\end{thebibliography}
\end{document}